\documentclass[12pt]{article}
\usepackage{mathrsfs}
\usepackage{amsfonts}
\usepackage{amsmath}
\usepackage{amssymb,amsmath}
\def\Z{\mathbb{Z}}

\def\C{\mathbb{C}}
\def\cl{\centerline}
\def\vs{\vspace*}
\def\H{\textbf{H}}
\def\sig{\sigma}
\def\Der{\mbox{\rm Der}}
\def\Inn{\mbox{\rm ad\,}}
\def\Aut{\mbox{\rm Aut\,}}
\def\IAut{\Gamma}
\def\QED{\hfill$\Box$}
\def\CL{\widehat{\mathcal L}}

\numberwithin{equation}{section}
\newtheorem{theo}{Theorem}[section]

\newtheorem{exam}[theo]{Example}
\newtheorem{coro}[theo]{Corollary}
\newtheorem{lemm}[theo]{Lemma}
\newtheorem{prop}[theo]{Proposition}

\newtheorem{rema}[theo]{Remark}

\begin{document}

\cl{{\large\bf Derivations and automorphism groups of }}
\cl{{\large\bf completed Witt Lie algebra}}\vs{6pt}

\cl{ Yongping Wu$^{\,*}$, Ying Xu$^{\,\dag}$, Lamei Yuan$^{\,\dag}$}
 \cl{\small$^{*}$School of Mathematics and Computer
Science, Longyan University, Fujian 364000, China} \cl{\small
$^{\dag}$Department of Mathematics, University of Science \!and
\!Technology \!of \!China, Hefei 230026, China}

\cl{\small E-mail: wyp\_504@sohu.com, xying@mail.ustc.edu.cn,
lmyuan@mail.ustc.edu.cn }\vs{6pt}

{\small\parskip .005 truein \baselineskip 3pt \lineskip 3pt

\noindent{{\bf Abstract.} In this paper, a simple Lie algebra,
referred to as  the completed
Witt Lie algebra, %$\CL$
is introduced. Its derivation algebra and automorphism group are
completely described.
%It is proved that all the derivations of $\CL$ are inner derivations and that all the
%automorphisms are inner automorphisms.
As a by-product, it is obtained that the first cohomology group of
this Lie algebra %$\CL$
with coefficients in its adjoint module is trivial. Furthermore, we
completely determine the conjugate classes of this Lie algebra under
its automorphism group, and also obtain that this Lie algebra does
not contain any nonzero ${\rm ad}$-locally finite element.
%, namely $\H^1(\CL,\CL)=0$.
\vs{5pt}

\noindent{\bf Key words:} Completed Witt Lie algebra, Derivations,
Automorphism groups.}

\parskip .001 truein\baselineskip 6pt \lineskip 6pt

\vs{10pt}

\cl{\bf\S1. \
Introduction}\setcounter{section}{1}\setcounter{equation}{0}
\vs{6pt}

The complex Witt algebra $W_1$, named after Ernst Witt, was first
defined by E. Cartan \cite{C} in 1909, which is the Lie algebra of
meromorphic vector fields defined on the Riemann sphere that are
holomorphic except at two fixed points. It is also the
complexification of the Lie algebra of polynomial vector fields on a
circle, and the Lie algebra of derivations of the Laurent polynomial
ring $\C[t,t^{-1}]$. Witt algebras occur in the study of conformal
field theory. There are some related Lie algebras defined over
finite fields, studied by Witt in the 1930s, that are also called
Witt algebras.

The complex Witt algebras were the first examples of nonclassical
infinite dimensional simple Lie algebra in 1937. It has been paying
more and more attentions since it was found. As is known that the
$1$-dimensional central extension of the Witt algebra $W_1$ is the
well-known Virasoro algebra, which plays important roles both in
mathematics and physics. The classification of Harish-Chandra
modules over the Virasoro algebras, high rank and higher rank
Virasoro algebras are well developed in \cite{O,Su2,Su1,Su3}.\par
 The simple modules of the Witt algebra in prime characteristic
 were classified by Chang \cite{HJ} in 1941. For a long time this was the only systematic knowledge of
the representations of nonclassical simple Lie algebras. In the mid
1980s, Shen  \cite{GY1} obtained the ``mixed product'' realization
over tensor product of vector spaces. In this way, a complete
description of their graded and filtered simple representations was
established in \cite{GY1,GY2,GY3}. A simple description of the
representations of the Witt algebra was given in \cite{HS}.

Generalized Witt algebras were defined by Kaplansky \cite{I} in the
context of the classification problem of simple finite dimensional
Lie algebras over fields of prime characteristic. Later on,
generalizations of the simple Lie algebras of Witt type over a field
of characteristic zero have been studied in \cite{N,JM,DS,X}. The
structure spaces and some representations of these simple Lie
algebras of generalized Witt type were presented in \cite
{SXZ,SZHOU}. The derivations, automorphisms and the second
cohomology groups of such kind of algebras have been studied by
several authors as indicated in the references (see e.g. \cite{ZZ,
SZ}). The bialgebras of generalized Witt type were also considered
in \cite{SS,WSS}.

The classical one-sided Witt algebra $W_1^{+}$ is defined with basis
$\{L_i\,|\,i\in\Z,\,i\ge-1\}$ and brackets
$$[L_{i},L_{j}]=(j-i)L_{i+j}  \ \ \mbox{for all}\ i,j\geq -1, \ i,j\in\Z. $$
Thus by realizing $L_i$ as $t^{i+1}\frac{d}{dt}$, one immediately
observes that $W_1^+=\Der\,\mathbb{C}[t]=\C[t]\frac{d}{dt}$ is the
derivation algebra of the polynomial algebra $\mathbb{C}[t]$. From
this, one sees that it is very natural to consider the
generalization of the Witt algebra by replacing the polynomial
algebra $\mathbb{C}[t]$ by the power series algebra
$\mathbb{C}[[t]]$. In this paper, by starting from the determination
of the derivation algebra of $\mathbb{C}[[t]]$, we are able to
define a Lie algebra $\CL$}, which we refer to as the {\it completed
Witt algebra}, as follows: Let ${L_{i}, i\geq-1,i \in \Z}$ be
symbols, denoted by $\CL$ the $\C$-vector space consisting of
elements $\sum_{i=-1}^{\infty}a_iL_i$ with $a_i\in\mathbb{C}$ such
that $\sum_{i=-1}^{\infty}a_iL_i=\sum_{i=-1}^{\infty}b_iL_i$ if and
only if $a_i=b_i$ for all $i\geq-1$. Then $\CL$ is a Lie algebra
under the brackets
\begin{eqnarray}\label{11}
\Big[\mbox{$\sum\limits_{i=-1}^{\infty}$}a_iL_i,\mbox{$\sum\limits_{i=-1}^{\infty}$}b_iL_i\Big]=
\mbox{$\sum\limits_{i=-1}^{\infty}$}\Big(\mbox{$\sum\limits_{j=-1}^{i}$}(2j-i)a_{i-j}b_{j}\Big)L_i.
\end{eqnarray}
One will see from \eqref{DERCC} that $\CL$ is in fact the derivation
algebra of $\C[[t]]$ and it is a simple Lie algebra (Theorem
\ref{Mai-1}).
 One of the motivations in studying this Lie algebra
also comes from a recent paper \cite{YS}, where in order to
determine a new quantization of the Witt algebra, it is necessary to
consider its completion. Our first aim is to study the structure
theory of this Lie algebra, such as the determinations of the
derivation algebra and the automorphism group of this Lie algebra.
The main results will be summarized in Theorems \ref{Main-Theo1},
\ref{Main-Theo2} and \ref{Main-Theo3}. As a by-product, we
completely determine the conjugate classes of $\CL$ under the
automorphism group, and also obtain that $\CL$ does not contain any
nonzero ${\rm ad}$-locally finite element (Corollary
\ref{Main-co3}). Furthermore, we shall also present some interesting
examples of subalgebras of this Lie algebra in Section 3; in
particular, the centerless Virasoro algebra and nongraded Virasoro
algebra, Lie algebras of Block type can be all realized as
subalgebras of this Lie algebra (Examples \ref{exam3}, \ref{exam4}
and \ref{exam6}). Our next aim will be the study of the
representations of this Lie algebra.

Finally, we would like to point out that since the algebras we are
considering are infinite-dimensional with uncountable basis (and in
fact it is impossible to find a basis of such an algebra), some
techniques developed in \cite{F} cannot be applied to our cases
here. Also,  since $\CL$ does not contain any nonzero locally finite
elements, some standard methods in determining automorphisms (such
as those used in \cite{SZhao}) cannot be applied either. Thus it
seems to us that it is necessary to develop some new techniques in
dealing with the problems occurring in the present paper. This is
also our another motivation to present this paper. \vs{12pt}\

\cl{\bf\S2. \ Derivation algebra of
$\mathbb{C}[[t]]$}\setcounter{section}{2}\setcounter{equation}{0}\setcounter{theo}{0}

\vs{8pt}

In this section, we shall determine the derivation algebra of the
Laurent polynomial ring
$\mathbb{C}[[t]]=\{\mbox{$\sum_{i=0}^{\infty}$}a_it^i\mid
a_i\in\mathbb{C}\}.$ Let us begin with the definition of
derivations.\par Let $\mathcal {A}$ be an algebra over $\mathbb{C}$.
A {\it derivation} $D$ of $\mathcal {A}$ is a linear transformation
on $\mathcal {A}$ such that
\begin{equation}\label{Der-1}D(xy)=D(x)y+xD(y)\ \ \mbox {for}\ \ x,y\in
\mathcal {A}.\end{equation}
 The space $\Der\,\mathcal {A}$ of
derivations forms a Lie algebra with respect to the Lie
bracket:$$[D_1,D_2]=D_1D_2-D_2D_1\ \ \mbox {for} \ \ D_1,D_2\in
\Der\, \mathcal {A}.$$ For any $p(t)\in\C[[t]]$, one can define a
derivation $D=p(t)\frac{d}{dt}\in\Der\,\C[[t]]$ as follows:
$$ D(a(t))=p(t)\frac{d}{dt}a(t)\in\C[t]] \mbox{ \ for any
\ }a(t)\in\C[[t]].$$ We denote
$\C[[t]]\frac{d}{dt}=\{p(t)\frac{d}{dt}\,|\,p(t)\in\C[[t]]\}$. The
main result in this section is the following.
\begin{theo}\label{Main-Theo1} The derivation algebra of $\C[[t]]$
is $\Der\,\C[[t]]=\C[[t]]\frac{d}{dt}$.
\end{theo}
\noindent{\it Proof.~}~Let $D\in\Der\,\C[[t]]$. We shall prove that
after a number of steps in each of which $D$ is replaced by $D-D'$
for some $D'\in \C[[t]]\frac{d}{dt}$, the 0 derivation is obtained
and thus proving that $D\in\C[[t]]\frac{d}{dt}$. This will be done
by two lemmas below.\hfill$\Box$
\begin{lemm}
Let $D\in \Der\, \mathbb{C}[[t]]$. By replacing $D$ by $D-D'$ for
some $D'\in\C[[t]]\frac{d}{dt}$, we can suppose $D(t^i)=0$ for
$i\ge0.$
\end{lemm}
\noindent{\it Proof}.\ \ Let  $D\in \Der\, \mathbb{C}[[t]]$. Take
$p(t)=D(t)\in \mathbb{C}[[t]]$. Re-denoting $D-p(t)\frac{d}{dt}$ by
$D$, we have $D(t)=0$. Using \eqref{Der-1} and induction on $i>0$,
we have $D(t^i)=0$. Obviously, we have $D(t^0)=D(t^0
t^0)=D(t^0)+D(t^0)$, and so $D(t^0)=0$.\hfill$\Box$
\begin{lemm}
Suppose $D\in\Der\,\C[[t]]$ satisfying $D(t^i)=0$ for $i\ge0$. Then
$D=0$.
\end{lemm}
\noindent{\it Proof}.\ \ Let $q(t)=\sum_{i=0}^\infty
a_it^i\in\C[[t]]$. Assume
$D(q(t))=\mbox{$\sum_{i=0}^{\infty}$}b_it^i$. To compute $b_{i_0}$
for any fixed $i_0\ge0$, we take a fixed $N>i_0$, and we have
\begin{eqnarray}\label{DEEE0}
D(q(t))&=&D\Big(\mbox{$\sum\limits_{i=0}^{N}a_{i}t^{i}$}
+\mbox{$\sum\limits_{i=N}^{\infty}a_{i}t^{i}$}\Big)=D\Big(\mbox{$\sum\limits_{i=0}^{N}a_{i}t^{i}$}\Big)
+D\Big(\mbox{$\sum\limits_{i=N}^{\infty}a_{i}t^{i}$}\Big)
\nonumber\\&=&
D\Big(\mbox{$\sum\limits_{i=N}^{\infty}a_{i}t^{i}$}\Big)=
D\Big(t^{N}\mbox{$\sum\limits_{i=0}^{\infty}a_{i+N}t^{i}\Big)$}=
t^{N}D\Big(\mbox{$\sum\limits_{i=0}^{\infty}a_{i+N}t^{i}\Big)$},
\end{eqnarray}
where the third and last quialities follow from \eqref{Der-1} and
the fact that $D(t^i)=0$ for $i\in\Z$. Note that the right-hand side
of \eqref{DEEE0}, being a power series, does not contain the term
$t^{i_0}$ since $N>i_0$. Thus $b_{i_0}=0$. Since $i_0$ is
arbitrarily chosen, we obtain $D(q(t))=0.$\QED\vs{7pt}\par

If we realize $L_i$ as $t^{i+1}\frac{d}{dt}$, then the Lie algebra
$\CL$ defined in \eqref{11} is simply
\begin{equation}\label{DERCC}\CL=\Der\,\mathbb{C}[[t]]=\mathbb{C}[[t]]\frac{d}{dt}\vs{-6pt}.\end{equation}
\begin{theo}\label{Mai-1}The Lie algebra $\CL$ is simple.\end{theo}
\noindent{\it Proof.~}~Let $I$ be any nonzero ideal of $\CL$, and
let $x=\sum_{i=i_0}^\infty a_{i}L_i\in I$ be a nonzero element,
where $i_0\ge-1$ is the smallest integer such that $a_{i_0}\ne0$.
Applying ${\rm ad\,}L_{-1}$ to $x$ several times and re-denoting the
result by $x$ if necessary, we can suppose $i_0=-1$. By rescaling
$x$, we can suppose $a_{-1}=1$. Now let $y=\sum_{i=-1}^\infty
b_iL_i\in\CL$ be any element. We can always find some
$z=\sum_{i=-1}^\infty c_iL_i\in\CL$ for some $c_i\in\C$ such that
$y=[x,z]$, which is equivalent to
\begin{equation}\label{oais}
b_i=\mbox{$\sum\limits_{j=-1}^{i+1}$}
(2j-i)a_{i-j}c_{j}=c_{i+1}+\mbox{$\sum\limits_{j=-1}^{i}$}(2j-i)a_{i-j}c_{j}\mbox{
\ for \ }i\ge-1.
\end{equation}
Regarding \eqref{oais} as a system of linear equation on
$c_{i},\,i\ge-1$, we see that there exists a unique solution for
$c_i,\,i\ge-1$, such that $y=[x,z]$ holds. This implies $y\in I$,
i.e., $I=\CL$.\QED \vskip7pt

 \cl{\bf\S3. \ Some interesting examples of subalgebras of
$\CL$}
\setcounter{section}{3}\setcounter{equation}{0}\setcounter{theo}{0}
\vs{8pt}

We present some interesting examples of subalgebras of $\CL$ below.
One may expect that their structure and representation theories are
of interests and worth further studies.
\begin{exam}\label{exam1}\rm \begin{itemize}\item[\rm(1)]
 Let $B_1={\rm Span}\{L_n=2\sin nt\frac{d}{dt}\,|\,n\in \Z\}$. Then
$L_i=-L_{-i}$ and
\begin{equation}\label{Ex--1}[L_i,L_j]=(j-i)L_{i+j}+(i+j)L_{i-j}.\end{equation}
\item[\rm(2)]
%\end{exam}\begin{exam}\label{exam2}\rm
Let $B_2={\rm Span}\{L_n=2\sin nt\frac{d}{dt},M_n=2\cos
nt\frac{d}{dt}\,|\,n\in \Z\}$. Then $L_i=-L_{-i}$, $M_i=M_{-i}$, and
we have \eqref{Ex--1} and
$$[L_i,M_j]=(j-i)M_{i+j}-(i+j)M_{i-j},
\ \ %$$  $$
[M_i,M_j]=(i-j)L_{i+j}+(i+j)L_{i-j}.$$
\end{itemize}\end{exam}\begin{exam}\label{exam3}\rm Let $B_3={\rm
Span}\{L_n=\exp nt\frac{d}{dt}\,|\,n\in \Z\}$. Then
$[L_i,L_j]=(j-i)L_{i+j}.$ Thus $B_3$ is the Witt algebra $W_1$ (or
the centerless Virasoro algebra).
\end{exam}\begin{exam}\label{exam4}\rm
 Let $B_4={\rm Span}\{L_{m,n}=\sin mt\exp
nt\frac{d}{dt},M_{m,n}=\cos mt\exp nt\frac{d}{dt}\,|\,m,n\in \Z\}$.
Then $L_{m,n}=-L_{-m,n},$ $M_{m,n}=M_{-m,n}$, and
$$[L_{m,n},L_{k,l}]=\frac{k-m}{2}L_{m+k,n+l}+\frac{k+m}{2}L_{m-k,n+l}+\frac{n-l}{2}M_{m+k,n+l}+\frac{l-n}{2}M_{m-k,n+l},$$
$$[L_{m,n},M_{k,l}]=\frac{k-m}{2}M_{m+k,n+l}-\frac{k+m}{2}M_{m-k,n+l}+\frac{l-n}{2}L_{m+k,n+l}+\frac{l-n}{2}L_{m-k,n+l},$$
$$[M_{m,n},M_{k,l}]=\frac{m-k}{2}L_{m+k,n+l}+\frac{k+m}{2}L_{m-k,n+l}+\frac{l-n}{2}M_{m+k,n+l}+\frac{l-n}{2}M_{m-k,n+l}.$$
\end{exam}\begin{exam}\label{exam5}\rm
Let $B_5={\rm Span}\{L_{m,n}=t^n\exp mt\frac{d}{dt}\,|\,m,n\in
\Z\}$. Then
$$[L_{m,n},L_{k,l}]=(k-m)L_{m+k,n+l}+(l-n)L_{m+k,n+l-1}.$$
Thus $B_5$ is in fact the centerless nongraded Virasoro defined in
\cite[(4.5)]{SZ}.
\end{exam}\begin{exam}\label{exam6}\rm \begin{itemize}\item[\rm(1)]
Let $B_6={\rm Span}\{L_{m,n}=t^m(t+a)^n\frac{d}{dt}\,|\,m,n\in \Z\}$
($a\neq 0$ is fixed). Then
$$[L_{m,n},L_{k,l}]=(k-m)L_{m+k-1,n+l}+(l-n)L_{m+k,n+l-1}.$$
\item[\rm(2)]
%\end{exam}\begin{exam}\label{exam7}\rm
Let $B_7={\rm
Span}\{L_{m,n}=(\sin t)^m(\cos t)^n \frac{d}{dt}\,|\,m,n\in \Z\}$.
Then
$$[L_{m,n},L_{k,l}]=(k-m)L_{m+k-1,n+l+1}-(l-n)L_{m+k+1,n+l-1}.$$
\end{itemize}
\noindent Thus $B_6,\,B_7$ are in fact Lie algebras of Block type
studied in \cite{CS-B,Li-S-B,Su-B1,Su-B2,WS-Bl,YueS}, etc.
\end{exam}\vskip7pt

 \cl{\bf\S4. \ Derivations of
$\CL$}\setcounter{section}{4}\setcounter{equation}{0}\setcounter{theo}{0}

\vskip8pt

In this section, we aim to determine all the derivations of the
completed Witt Lie algebra $\CL$  with the Lie bracket defined by
equation (\ref{11}). As in \eqref{Der-1}, a linear transformation
$D:\,\CL\rightarrow \CL$ is a
 {\it derivation} of $\CL$, if it satisfies  $D( [x,y])=[D
(x),y]+[x,D(y)]$ for $x,y\in\CL.$ If, in addition, there exists some
$x\in\CL$ such that $D(y)={\rm ad\,}x(y)=[x,y]$ for $y\in\CL,$ then
$D={\rm ad\,}x$ is called an inner derivation. We denote by
$\Der\,\CL$ (resp., $\Inn\CL$) the set of all derivations (resp.,
inner derivations). Then the main result in this section is the
following.
\begin{theo}\label{Main-Theo2} The derivation algebra of $\CL$
is $\Der\,\CL=\Inn\CL$. In particular, the first cohomology group of
$\CL$ with coefficients in its adjoint module is trivial, namely,
$\H^{1}(\CL,\CL)=0$.
\end{theo}
\noindent{\it Proof.~}~Let $D\in\Der\,\CL$. Again we shall prove
that after a number of steps in each of which $D$ is replaced by
$D-D'$ for some $D'\in \Inn\CL$, the 0 derivation is obtained and
thus proving that $D\in\Inn\CL$. This will be done by three lemmas
below.\hfill$\Box$
\begin{lemm}\label{lmma1}
Let $D\in \Der\, \CL$, by replacing $D$ by $D-D'$ for some
$D'\in\Inn\CL$, we can suppose  $D(L_{i}) = a_{i}L_{i}$ for some
$a_i\in\C$ with $i\geq-1$ and $a_0=0$.
\end{lemm}

\noindent{\it Proof}.\ \ Let $D \in \Der\,\CL$. Denote
$D(L_{0})=\sum_{i=-1}^{\infty}a_{i}L_{i}.$ Taking
$v=\sum_{i\neq0}(-\frac{a_{i}}{i})L_{i}\in\CL$, one has
$$
(D-{\rm ad\,} v)(L_{0})=D(L_{0})-[v,L_{0}]
=\mbox{$\sum\limits_{i=-1}^{\infty}$}a_{i}L_{i}-
\mbox{$\sum\limits_{i\neq0}$}(-\frac{a_{i}}{i})[L_{i},L_{0}]=a_{0}L_{0}.
$$
Re-denoting $D-{\rm ad\,}v$ by $D$, we have $D(L_{0})=a_{0}L_{0}$.
Suppose that
$D(L_{i})=\mbox{$\sum\limits_{j=-1}^{\infty}a_{ij}L_{j}$}\
(i\neq0)$. Applying $D$\ to $[L_{0},L_{i}]=iL_{i}$,\ we have
$ia_{0}L_{i}+\mbox{$\sum_{j\neq i}a_{ij}jL_{j}$}=i\mbox{$\sum_{j\neq
i}a_{ij}L_{j}$},$ which implies $a_{0}=0$ since $i\neq0$, and
$a_{ij}=0$\ for $j\neq i$. Namely, $D(L_{0})=0$\ and
$D(L_{i})=a_{ii}L_{i}, \,i\neq 0$. Let $a_i=a_{ii}$ for  $i\geq -1$,
then we get $D(L_{i})=a_{i}L_{i}$ with $a_0=0.$
 \QED
\begin{lemm}\label{lmma2}
Let $D\in \Der\, \CL$, by replacing $D$ by $D-D'$ for some
$D'\in\Inn\CL$, we can suppose  $D(L_{i}) = 0$ for all $i\geq-1$.
\end{lemm}
\noindent{\it Proof}.\ \ Suppose $D\in\Der\,\CL$ satisfies Lemma
\ref{lmma1}. Re-denoting $D-a_{1}{\rm ad\,} L_{0}$ by $D$, we obtain
$D(L_0)=D(L_1)=0$. Applying $D$ to $[L_{i},L_{1}]=(1-i)L_{i+1}$, we
obtain $a_{i}(1-i)L_{i+1}=(1-i)a_{i+1}L_{i+1},$ which implies
\begin{eqnarray}\label{21}
a_{i}=a_{i+1} \ \ \mbox{if} \ \ i\neq1.
\end{eqnarray}
So $a_{-1}=a_{0}=0$. Applying $D$ to $[L_{-1},L_{2}]=3L_{1}$, we get
$[L_{-1},a_{2}L_{2}]=0$, i.e. $3a_{2}L_{1}=0$, and thus $a_{2}=0$.
Then the result follows from (\ref{21}).\QED
\begin{lemm}\label{lmma3}
Suppose $D\in\Der\,\CL$ satisfying $D(L_i)=0$ for $i\ge-1$. Then
$D=0$.
\end{lemm}
\noindent{\it Proof}.\ \ Let
$x=\sum_{i=-1}^{\infty}a_{i}L_{i}\in\CL$. Assume
$D(x)=\sum_{i=-1}^{\infty}b_{i}L_{i}$.\ To compute $b_{i_0}$ for a
fixed $i_0\ge-1$, we take a fixed $N>i_0+1$, and we have
\begin{eqnarray*}\label{MAMAMA}
D(x)&=&D\Big(\mbox{$\sum\limits_{i=-1}^{\infty}a_{i}L_{i}$}\Big)
=D\Big(\mbox{$\sum\limits_{i=-1}^{N-2}a_{i}L_{i}$}
+\mbox{$\sum\limits_{i=N-1}^{\infty}a_{i}L_{i}$}\Big)=
D\Big(\mbox{$\sum\limits_{i=N-1}^{\infty}a_{i}L_{i}$}\Big)
\\&=&D\Big(\mbox{$\sum\limits_{i=-1}^{\infty}a_{i+N}L_{i+N}$}\Big)=
D\Big(\Big[L_{N},\mbox{$\sum\limits_{i\neq
N}\frac{a_{i+N}}{i-N}L_{i}\Big]+a_{2N}L_{2N}\Big)$}\\&=&
\Big[L_{N},D\Big(\mbox{$\sum\limits_{i\neq
N}\frac{a_{i+N}}{i-N}L_{i}\Big)$}\Big].
\end{eqnarray*}
If we write $D(\mbox{$\sum_{i\neq
N}\frac{a_{i+N}}{i-N}L_{i})$}=$\mbox{$\sum_{i=-1}^{\infty}c_{i}L_{i}$},
then we obtain
$\mbox{$\sum_{i=-1}^{\infty}b_{i}L_{i}$}=[L_{N},\mbox{$\sum_{i=-1}^{\infty}c_{i}L_{i}$}]
=\mbox{$\sum_{i=-1}^{\infty}c_{i}(i-N)L_{i+N}$},$ which in
particular implies $b_{i_0}=0$.  Since $i_0$ is arbitrarily chosen,
we obtain  $D(x)=0$.\QED\vskip 8pt

 \cl{\bf\S5. \ Automorphism Groups of
$\CL$}\setcounter{section}{5}\setcounter{equation}{0}\setcounter{theo}{0}

\vs{8pt}

In this section we compute the automorphism group  $\Aut\CL$ of the
completed Witt algebra $\CL$. This will be done by several lemmas.
We remark that since $\CL$ does not contain any nonzero locally
finite elements (Corollary \ref{Main-co3}), some standard methods
such as those used in \cite{SZhao} cannot be applied in our case
here.
\begin{lemm}\label{lmma5} For any $\sig\in \Aut\CL$, write
$\sig(L_{0})=\mbox{$\sum_{i=-1}^{\infty} $}a_{i}L_{i}$. Then one has
$(a_{-1},a_{0})\neq(0,0)$.
\end{lemm}
\noindent{\it Proof}.\ \ If not the case, then
$\sig(L_{0})=\mbox{$\sum_{i=1}^{\infty} $}a_{i}L_{i}$. Write
$\sig(L_{-1})=\mbox{$\sum_{j=j_{0}}^{\infty} $}b_{j}L_{j}$,\ where
$j_{0}$ is the smallest index such that $b_{j_{0}}\neq 0$. Applying
$\sig$ to $[L_{0},L_{-1}]=-L_{-1}$, one has
$[\sig(L_{0}),\sig(L_{-1})]=-\sig(L_{-1})$,
namely,{\def\Big{}\def\limits{}
$\Big[\mbox{$\sum\limits_{i=1}^{\infty}
$}a_{i}L_{i},\mbox{$\sum\limits_{j=j_{0}}^{\infty}
$}b_{j}L_{j}\Big]=-\mbox{$\sum\limits_{j=j_{0}}^{\infty}
$}b_{j}L_{j}.$} Note that the left-hand side does not contain the
term $L_{j_0}$, which implies that $b_{j_{0}}=0$, a contradiction
with the assumption that $b_{j_0}\ne0$.\QED\vskip7pt\par For
convenience, we introduce the following {\it filtration} of $\CL$:
\begin{eqnarray}\label{fil}
\CL=\CL_{-1}\supset\CL_{0}\supset\CL_{1}\supset\cdots\mbox{ \ with \
}\CL_{i}=\Big\{\mbox{$\sum\limits_{j=i}^{\infty} $}a_{j}L_{j}\mid
a_{j}\in \C\Big\},\end{eqnarray} which satisfies
\begin{eqnarray}\label{fil-2}
\bigcap_{i=-1}^{\infty}\CL_{i}=0, \ \ \
[\CL_i,\CL_j]\subset\CL_{i+j}\mbox{ \ for \ }i,j\ge-1.
%$\bigcup_{i=-1}^{\infty}\CL_{i}=\CL$
\end{eqnarray}

\begin{rema}\label{rema-1}\rm
For any $x=\sum_{i=-1}^{\infty} a_{i}L_{i}\in\CL$, if we need to
determine the coefficient $a_{i_0}$ of $L_{i_0}$ in $x$ for a fixed
$i_0$, we can always do the computation under modulo $\CL_{N}$ for a
fixed $N\gg0$. In this way, we can simplify the computation by
reducing a sum with infinite terms to a sum with finite terms.
\end{rema}

For any $x=\sum_{i=1}^{\infty} b_{i}L_{i}\in\CL_1$, we can define a
linear transformation $\exp^{{\rm ad\,}x}$ on $\CL$ by
\begin{equation}\label{exp-1} \exp^{{\rm
ad\,}x}(y)=\sum_{k=0}^\infty\frac1{k!}({\rm ad\,}x)^k(y)\mbox{ \ for
\ }y\in\CL.\end{equation} Note that for any $y\in\CL$, the
right-hand side of \eqref{exp-1} well defines an element  in $\CL$
by observing the following: for any fixed $N$, since $x\in\CL_1$, by
\eqref{fil-2}, we have $({\rm ad\,}x)^k(y)\in\CL_N$ when $k>N$, thus
Remark \ref{rema-1} can be applied.
\begin{lemm}\label{lmma6}
If $x\in \CL_1$, then $\exp^{{\rm ad\,}x}\in \Aut\CL$.
\end{lemm}
\noindent{\it Proof}.\ \ Obviously, $\exp^{{\rm ad\,}x}$ is a
bijection, since one can immediately check that $\exp^{{\rm
ad\,}(-x)}\exp^{{\rm ad\,}x}$ is the identity transformation. Thus
it remains to prove \begin{equation}\label{exp-2} [\exp^{{\rm
ad}\,x}(y),\exp^{{\rm ad\,}x}(z)]=\exp^{{\rm ad\,}x}([y,z])\mbox{ \
for \ }y,z\in\CL.\end{equation} By \eqref{fil-2} and Remark
\ref{rema-1}, we can suppose $x,y,z$ and the right-hand side of
\eqref{exp-1} all contain only finite terms. In this case,
\eqref{exp-2} is obvious. \QED\vskip7pt\par Similarly, we can prove
that $\exp^{{\rm ad\,}L_0}\in\Aut\CL.$ Denote
\begin{eqnarray}\label{iaut}
\IAut=\mbox{ \ the subgroup of $\Aut\CL$ generated by
$\exp^{a_0\,{\rm ad\,}L_{0}}$ and $\exp^{{\rm ad\,}x}$ for
$a_0\in\C,\ x\in\CL_1$.}
\end{eqnarray}
\begin{prop}\label{Main-Theo3-} The automorphism group of $\CL$ is $\Aut\CL=\IAut$.
\end{prop}
\noindent{\it Proof.~}~The proof will be done by several lemmas
below.\hfill$\Box$\vskip7pt
\begin{lemm}\label{lmma7} With notions as in Lemma $\ref{lmma5},$
there must exist $x\in\CL_1$ such that $\sig (L_{0})=\exp^{{\rm
ad\,}x}(a_{0}L_{0}+a_{-1}L_{-1})$. Thus by re-denoting $\exp^{-{\rm
ad\,}x}\sig $ by  $\sig$, we can now suppose
$\sig(L_0)=a_{-1}L_{-1}+a_0L_0$.
\end{lemm}
\noindent{\it Proof}.\ \ We need to prove that there exists
$x=\mbox{$\sum_{i=1}^{\infty}$}b_{i}L_{i}\in\CL_1$ for some
$b_i\in\C$  such
that\begin{eqnarray}\label{Solve}\!\!\!\!\!\!\!\!\!\!\mbox{$
\sum\limits_{i=-1}^{\infty}$}a_{i}L_{i}\!\!\!\!&=\!\!\!\!&\sig
(L_{0})=\exp^{{\rm
ad\,}x}(a_{0}L_{0}+a_{-1}L_{-1})=\mbox{$\sum\limits_{i=0}^{\infty}$}\frac{1}{i!}({\rm ad\,}x)^{i}(a_{0}L_{0}+a_{-1}L_{-1})\nonumber\\
&=\!\!\!\!&(a_{0}L_{0}+a_{-1}L_{-1})+[x,a_{0}L_{0}+a_{-1}L_{-1}]
+\frac{1}{2!}[x,[x,[a_{0}L_{0}+a_{-1}L_{-1}]]]+\cdots\end{eqnarray}
By comparing the coefficients of $L_j$ for $j\geq-1$, we obtain the
following system of equations on unknown variables $b_i,\,i\ge1$:
\begin{eqnarray}\label{Solve1}
\left\{
\begin{array}{lll} a_0&=&a_0-2a_{-1}b_1,\\
a_1&=&-a_0b_1-3a_{-1}b_2+a_{-1}b_1^2,\\
a_2&=&-2a_0b_2-4a_{-1}b_3+2a_{-1}b_1b_2,\\
a_3&=&-3a_0b_3-5a_{-1}b_4+a_{-1}b_1b_3-\frac{1}{2}a_0b_1b_2+\frac{3}{2}a_{-1}b_2^2+\frac{1}{3}a_{-1}b_1^2b_2,\\
 \cdots
\end{array}
\right.
\end{eqnarray}
We see that no matter whether $a_{-1}\ne0$ or $a_{-1}=0,\,a_0\ne0$,
the system \eqref{Solve1} has a unique solution for $b_i,\,i\ge1$,
which can be expressed in terms of $a_j$ with $-1\leq j<i$. Hence we
can choose $b_i,\,i\ge1$ such that \eqref{Solve} holds. The lemma
follows.\QED\vskip7pt\par

The following Lemma will play a key role in proving our main result.
\begin{lemm}\label{lmma9}
For any $x=\mbox{$\sum _{i=-1}^{\infty} $}a_{i}L_{i}\in \CL$,
suppose  $\sig(x)=\mbox{$\sum _{i=-1}^{\infty} $}a_{i}'L_{i}$. Then
for every $i_0$, the coefficient $a_i'$ of $L_{i_{0}}$ in $\sig(x)$
is the same as that in $\mbox{$\sum _{i=-1}^{N-1}
$}a_{i}\sig(L_{i})$ for $N\gg0$.
\end{lemm}
\noindent{\it Proof}.\ \ Let $i_0$ be fixed. Let
$L'_{i_0+2}:=\sig^{-1}(L_{i_0+2})=\sum_{j=j_0}^\infty b_{j}L_j$,
where $j_0\ge-1$ is the smallest integer such that $b_{j_0}\ne0$.
For any $N\ge{\rm max}\{i_0,2j_0+1\}$, we want to prove there exists
$y=\sum_{i=N-j_0}^\infty c_i L_i\in\CL$ for some $c_i\in\C$ such
that
\begin{equation}\label{0a0a-1}\mbox{$\sum\limits_{i=N}^\infty$}
a_iL_i=[y,L'_{i_0+2}],\end{equation}
 which is equivalent to
\begin{equation}\label{0a0a}
a_i=\mbox{$\sum\limits_{j=N-j_0}^{i-j_0}$}
(i-2j)c_jb_{i-j}=(2j_0-i)c_{i-j_0}b_{j_0}+\mbox{$\sum\limits_{j=N-j_0}^{i-j_0-1}$}(i-2j)c_jb_{i-j}\mbox{
\ for \ }i\ge N.
\end{equation}
Regarding \eqref{0a0a} as a system of linear equations on
$c_{i-j_0},\,i\ge N$, and noting that the coefficient
$(2j_0-i)b_{j_0}$ of $c_{i-j_0}$ is always nonzero since
$b_{j_0}\ne0$ and $i\ge N>2j_0$, we obtain that there exists q
unique solution of $c_i,\,i\ge N-j_0$ such that \eqref{0a0a-1}
holds.
\par Now using \eqref{0a0a-1}, we have
\begin{eqnarray*}
\sig(x)
&=&\sig\Big(\mbox{$\sum\limits_{i=-1}^{N-1}$}a_{i}L_{i}+\mbox{$\sum\limits_{i=N}^{\infty}
$}a_{i}L_{i}\Big)
=\mbox{$\sum\limits_{i=-1}^{N-1}$}a_{i}\sig(L_{i})+\sig\Big(\mbox{$\sum\limits_{i=N}^{\infty}
$}a_{i}L_{i}\Big)\\
&=&\mbox{$\sum\limits_{i=-1}^{N-1}$}a_{i}\sig(L_{i})+\sig([y,L'_{i_0+2}])=
\mbox{$\sum\limits_{i=-1}^{N-1}$}a_{i}\sig(L_{i})+[\sig(y),L_{i_0+2}].
\end{eqnarray*}
It is easy to see that the smallest index in the second term is
equal or greater than $i_0+1$. Hence, the coefficient $a'_{i_0}$ of
$L_{i_0}$ in $\sig(x)$ is the same as that in
$\mbox{$\sum_{i=-1}^{N-1} $}a_{i}\sig(L_{i})$.\QED\vskip7pt

\begin{lemm}\label{Claim1}We have $a_{-1}=0$ and $a_0=1.$ \end{lemm}\noindent{\it Proof.~}
Write $\sig(L_i)=\sum_{j=-1}^\infty b_{i,j}L_j$ for $i\ge-1$.
Applying $\sig$ to $iL_i=[L_0,L_i]$ and comparing the coefficients
of $L_j$ in both sides, we obtain
\begin{equation}\label{AJAJ}ib_{i,j}=(j+2)a_{-1}b_{i,j+1}+ja_0b_{i,j}\mbox{
\ for \ }i,j\ge-1.\end{equation}First assume $a_{-1}\ne0$. From
\eqref{AJAJ}, we inductively obtain \begin{equation}\label{AJAJ1}
b_{i,j}=\binom{\frac{i}{a_0} + 1}{j +
1}\Big(\frac{a_0}{a_{-1}}\Big)^{j+1} b_{i,-1}\mbox{ \ and \ }
b_{i,-1}\ne0\mbox{ \ for \ }i,j\ge-1,
\end{equation}
where in general, the binomial coefficient $\binom{a}{j}$ is defined
to be $\frac{a(a-1)\cdots(a-j+1)}{j!}\vs{2pt}$ if $0\le j\in\Z$, or
zero otherwise, and where, the second equation follows from the
first and the fact that $\sig(L_i)\ne0$. Note that \eqref{AJAJ1} has
meaning even if $a_0=0$ as we can rewrite \eqref{AJAJ1} so that
$a_0$ does not appear in the denominator. Applying $\sig$ to
$(j-i)L_{i+j}=[L_i,L_j]$, and comparing the coefficients of $L_p$ in
both sides, we obtain
\begin{equation}\label{AJAJ2} (j - i)b_{i + j,
p}=\mbox{$\sum\limits_{k=-1}^{p+1}$}(p - 2 k) b_{i, k} b_{j, p -
k}\mbox{ \ for \ }i,j,p\ge-1.
\end{equation}
Taking $p=-1$ and using \eqref{AJAJ1}, we obtain
\begin{equation}\label{AJAJ3}
a_{-1}b_{i + j, -1}=b_{i,-1}b_{j,-1}\mbox{ \ for \ $i,j\ge-1$ \ with
\ }i\ne j.
\end{equation}
Taking $j=-1$, we see that $b_{-1,-1}\ne0$. From this, one
immediately solves
\begin{equation}\label{AJAJ4-}
b_{i,-1}=\Big(\frac{a_{-1}}{b_{-1,-1}}\Big)^{i+1}b_{-1,-1}\mbox{ \
for  \ }i\ge-1.\end{equation}
%This together with \eqref{AJAJ4} gives
%\begin{equation}\label{AJAJ4}
%b_{i,j}=\binom{\frac{i}{a_0} + 1}{j +
%1}\Big(\frac{a_0}{a_{-1}}\Big)^{j+1}(\frac{a_{-1}}{b_{-1,-1}})^{i+1}b_{-1,-1}\mbox{
%\ for \ }i,j\ge-1.
%\end{equation}
Now take $x=\sum_{i=-1}^\infty \frac1{b_{i,-1}}L_i\in\CL$. Let us
 compute the coefficient, denoted by
$\lambda$, of $L_{-1}$ in $\sig(x)$. By Lemma \ref{lmma9}, $\lambda$
is equal to the coefficient of $L_{-1}$ in $\sum_{i=-1}^N
\frac1{b_{i,-1}}\sig(L_{i})$ for all $N\gg0$, i.e., $\lambda=N$ for
all $N\gg0$, which is impossible since $\lambda$ is a constant
number.

This proves $a_{-1}=0$. Thus $a_0\ne0$. Now \eqref{AJAJ} shows that
$b_{i,j}=0$ if $j\ne\frac{i}{a_0}$. However, for each $i_0\ge-1$,
there exists at least one integer $j_0\ge-1$ such that
$b_{i_0,j_0}\ne0$. This shows that $j_0$ must be $\frac{i_0}{a_0}$,
thus $\frac{i_0}{a_0}\ge-1$ is an integer for all $i_0\ge-1$. In
particular, $a_0=1$. \hfill$\Box$\vskip7pt

Denote $b_i=b_{i,i}$. Then \eqref{AJAJ} implies
\begin{equation}\label{MAMAM}\sig(L_i)=b_iL_i\mbox{ \ \ for \ \
}i\ge-1\mbox{ \ with \ }b_0=1.\end{equation} Applying $\sig$ to
$(j-i)L_{i+j}=[L_i,L_j]$, we can easily solve that $b_i=b^i$ for
$i\ge-1$, where $b=b_1$.

\begin{lemm}\label{lmma10} For any $\sig\in \Aut\CL$, by re-denoting
${\rm exp}^{a_0\,{\rm ad\,}L_0}\sig$ by $\sig$ for some $a_0\in\C$,
we can suppose $\sig(L_{i})=L_{i}$ for $i \geq -1$. Furthermore
$\sig=1$.
\end{lemm}
\noindent{\it Proof}.\ \ Taking $a_0=-\ln b\in\C$, one can easily
check
$$\exp^{a_0\,{\rm ad\,}L_0}\sig(L_1)=\exp^{a_0\,{\rm ad\,}L_0}(bL_1)=L_1.$$
Re-denoting  $\exp^{a_0\,{\rm ad\,}L_0}\sig$ by $\sig$, one can
assume that $\sig(L_{1})=L_{1}$, i.e., $b=1$. Thus,
$\sig(L_{i})=L_{i}$\ for  $i \geq -1$.
 Now using Lemma
\ref{lmma9}, we obtain that $\sig(x)=x$ for any
$x=\mbox{$\sum_{i=-1}^{\infty} $}a_{i}L_{i}\in \CL$, i.e., $\sig=1$.
This completes the proof of of Lemma \ref{lmma10}.\QED\vskip7pt From
Lemmas \ref{lmma7} and \ref{lmma10}, we in fact have the following.
\begin{theo}\label{Main-Theo3} Let
$\CL_1$ be defined as in \eqref{fil}. Then the automorphism group of
$\CL$ is $\Aut\CL=\{\exp^{{\rm ad\,}x}\exp^{a_0\,{\rm
ad\,}L_0}\,|\,a_0\in\C,x\in\CL_1\}$.
\end{theo}

 As the application of the above proof, we can obtain
Corollary \ref{Main-co3}. First we give a concepts.
 An element $s\in\CL$ is called
 {\it $\rm ad$-locally finite} if for any
given $v\in \CL$, the subspace ${\rm Span}\{({\rm ad\,} s)^m\cdot
v\,|\,m\in\Z_+\}$ of $\CL$ is finite-dimensional. Two elements $y$
and $z$ in $\CL$ are said to be {\it $(\Aut\CL)$-conjugate} or {\it
conjugate under the automorphism group of $\CL$} if $y\in(\Aut\CL)
z$.

\begin{coro}\label{Main-co3}
\begin{itemize}\item[\rm(1)]Two elements $y$
and $z$ in $\CL$ are conjugate under the automorphism group of $\CL$
if and only if there exists some $i\ge-1$ such that
$y,z\in\CL_i\backslash\CL_{i+1}$ and moreever in case $i=0$,
$y-z\in\CL_1$.
\item[\rm(2)] The Lie algebra $\CL$
does not contain any nonzero ${\rm ad}$-locally finite element.
\end{itemize}
\end{coro}
\noindent{\it Proof.~}~Let $y=\sum_{i=i_0}^\infty a_iL_i$ be any
nonzero element in $\CL$ such that $i_0\ge-1$ is the smallest
integer with $a_{i_0}\ne0$. As in the proof of Lemma \ref{lmma7}, we
can choose some $\sig=\exp^{{\rm ad\,}x}\in\Aut\CL$ with $x\in\CL_1$
such that $\sig(x)=a_{i_0}L_{i_0}$. Furthermore,
$a_{i_0}L_{i_0}=\sig_0(L_{i_0})$ for some $\sig_0\in\Aut\CL$ if and
only if $a_{i_0}=1$, or $i_0\ne1$ and $\sig_0=\exp^{\lambda\,{\rm
ad\,}L_0}$ for $\lambda=\frac{{\rm ln\,}a_{i_0}}{i_0}$. From this,
we obtain (1).

Now suppose $y$ is a nonzero ${\rm ad}$-locally finite element. Then
$L_{i_0}=\frac1{a_{i_0}}\sig(y)$ is also an ${\rm ad}$-locally
finite element. But for $z=\sum_{i=-1}^\infty L_i$, one can easily
check that $({\rm ad}\,L_{i_0})^k(z)=\sum_{j={\rm
max}\{k(i_0-1),-1\}}^\infty(j-2i_0)^kL_j$ for $k=0,1,...$, are
linear independent, which is a contradiction. This completes the
proof of the corollary.\QED \vskip12pt

\small\noindent{\bf Acknowledgements~}~This research was conducted
during the first author's visit to the University of Science and
Technology of China. She would like to thank the hospitality and the
financial support of the Department of Mathematics of USTC. Part of
this research was also supported by NSF grants 10825101 of China.
%\vs{8pt}

\end{document}